\documentclass{article}

\usepackage{arxiv}

\usepackage[utf8]{inputenc} % allow utf-8 input
\usepackage[T1]{fontenc}    % use 8-bit T1 fonts
\usepackage[colorlinks=true,urlcolor = black,citecolor=black,linkcolor=blue]{hyperref}
\usepackage{hyperref}
\usepackage{url}            % simple URL typesetting
\usepackage{booktabs}       % professional-quality tables
\usepackage{amsfonts}       % blackboard math symbols
\usepackage{nicefrac}       % compact symbols for 1/2, etc.
\usepackage{microtype}      % microtypography

\usepackage{lipsum}         % Can be removed after putting your text content
\usepackage{graphicx}
\usepackage{doi}

\usepackage[utf8]{inputenc}
\usepackage{amsmath}
\usepackage{amssymb}
\usepackage{color}
\usepackage{amssymb,amsfonts,amsmath,amsthm,amsbsy} % Math symbols.
\usepackage{tabularx} % Helpful for full text width tables.
\usepackage{mathtools,mathrsfs}
\usepackage{mhchem}
\usepackage{enumerate}
\usepackage{booktabs}
\usepackage{bm}
\usepackage{bbm}
\usepackage{multirow}
\usepackage{threeparttable}
\usepackage{float}

\usepackage[table,dvipsnames, svgnames]{xcolor}

\usepackage[stable]{footmisc}
\usepackage[sort&compress]{natbib}

\usepackage{indentfirst}
\usepackage{etoolbox,siunitx}

\usepackage{afterpage}
\robustify\bfseries
\usepackage{adjustbox}
\usepackage{chngpage}
\usepackage{rotating}
\usepackage{hvfloat}
\usepackage{soul}
\usepackage{soul}
\RequirePackage{color}

\usepackage{subcaption}
\usepackage{caption}

\usepackage{cleveref}       % smart cross-referencing

\usepackage{amsthm}   % For theorem-like environments
\newtheorem{theorem}{Theorem}
\newtheorem{proposition}{Proposition}

\title{On the Generalizations of the Cauchy-Schwarz-Bunyakovsky Inequality with Applications to Elasticity}

\author{ 
	Dimitra Labropoulou\\
	Department of Chemical Engineering\\ 
	University of Patras\\
	Patra 26504, Greece \\
	\texttt{dlabropoulou@chemeng.upatras.gr} \\
	\And
	Thanasis Labropoulos\\
	Department of Mathematics\\
	School of Applied Mathematics and Physical Sciences\\
	National Technical University of Athens\\
	Athens 15773, Greece\\	
	\texttt{tlabropoulos@mail.ntua.gr} \\
	\And
	Panayiotis Vafeas\thanks{Corresponding Author: vafeas@chemeng.upatras.gr}\\
	Department of Chemical Engineering\\ 
	University of Patras\\
	Patra 26504, Greece \\
	\texttt{vafeas@chemeng.upatras.gr} \\
	\And
\href{https://orcid.org/0000-0002-4788-9877}{\includegraphics[scale=0.06]{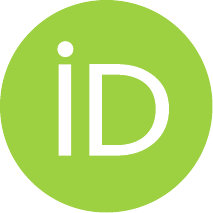}\hspace{1mm}Dimitris M. Manias}\\
	Department of Mathematics\\
	Khalifa University of Science and Technology\\
	Abu Dhabi 127788, United Arab Emirates \\
	\texttt{dimitris.manias@ku.ac.ae} \\
}

\renewcommand{\headeright}{A Postprint}
\renewcommand{\undertitle}{A Postprint}
\renewcommand{\shorttitle}{\textit{On the Generalizations of the Cauchy-Schwarz-Bunyakovsky Inequality } }

\hypersetup{
pdftitle={On the Generalizations of the Cauchy-Schwarz-Bunyakovsky Inequality with Applications to Elasticity},
pdfsubject={math.FA},
pdfauthor={Dimitris M. Manias},
}

\begin{document}
\maketitle

\begin{abstract}
In this article we present both the discrete and the integral form of Cauchy-Bunyakovsky-Schwarz ($\rm CBS$) inequality, some important generalizations in the $n$-dimensional Euclidean space and  in
linear subspaces of it, as well as the strengthened $\rm CBS$. The last $\rm CBS$ inequality plays an important role in elasticity problems. A geometrical interpretation and a collection of the most important proofs of it are, also, presented.
\end{abstract}

% keywords can be removed

\section{Introduction}
A.-L. Cauchy~\cite{Cau} in his monograph published in 1821,
presented an inequality (for finite discrete sums),  which
subsequently played and continues to play an important role in
mathematics. In 1859, V.I. Bunyakovsky~\cite{Bun} extended
this inequality to its integral version (that is, for the case
where we have continuous summation). Some years later, in 1888,
H.A. Schwarz~\cite{Sch} established the general form of
this inequality, valid in vector spaces endowed with dot product,
the so-called "Cauchy-Bunyakovsky-Schwarz inequality" (abbreviated
CBS inequality). 

The CBS inequality is one of the most famous
inequalities in mathematics, on the one hand, because of its close
relations with a multitude of other important inequalities of
mathematics and, on the other hand, because of its utility in a
particularly wide and varied range of applications in almost all
branches of classical and modern mathematics, such as Real and
Complex Analysis, Numerical Analysis, Hilbert Space Theory,
Differential Equations, Probability Theory, Statistics, Number
Theory, etc. More details on this topic we can see in the works of Dragomir\cite{Dra} and Steele\cite{Ste}. 

Additionally, new applications
of various branches of mathematics are constantly appearing, where
the CBS inequality itself or CBS-type inequalities play a significant
role. B${\rm \hat i}$rsan and B${\rm \hat i}$rsan~\cite{Bi-Bi} presented a CBS-type inequality in
both the discrete and the integral forms. The integral
version of this inequality appears in the study of mechanical 
properties of thin elastic rods. Also, Achchab and Maitre\cite{Ac-Ma} presented the Strengthened Cauchy - Bunyakowski - Schwarz inequality scheme, which plays a fundamental role in the convergence rate of multilevel iterative methods. Its optimal
constant plays an important role in elasticity problems (see
Margenov\cite{Mar}, Achchab et al.\cite{Ach-Axe}, Achchab and Maitre\cite{Ac-Ma}, Dragomir\cite{Dra}, Jung and Maitre\cite{Ju-Ma}). Besides, this handy theory can be readily applied to more recent developments in anisotropic elasticity, wherein even nowadays novel methods have been introduced in the work of Labropoulou et al.\cite{Lab}.

In this paper we attempt a presentation of the basic
versions of the CBS inequality, an inequality of CBS type, the
integral version of which appears in the study of mechanical
properties of thin elastic rods~\cite{Bi-Al, Bi-Bi}, demonstrating that the Strengthened CBS
inequality is very useful in elasticity problems. Finally,
in the Appendix we present an anthology of clever and fairly short
proofs of the CBS inequality (see also the work of Wu and Wu~\cite{Sh}).

\section{Cauchy-Bunyakovsky-Schwarz Inequality}

\subsection{Cauchy Inequality}
Cauchy's inequality for non-negative real numbers is expressed as
follows.
\begin{proposition} \label{pro2.1}\rm
If $\alpha _1 ,\,\alpha _2 ,\,...,\alpha _n$ are non-negative real numbers, then
\begin{equation}\label{eq1}
\frac{{\alpha _1  + \alpha _2  + ... + \alpha _\nu  }} {\nu }
\geqslant \sqrt[\nu ]{{\alpha _1  \cdot \alpha _2 ... \cdot \alpha
_\nu }}~.
\end{equation}

\noindent If in addition the numbers  $\alpha _1 ,\,\alpha _2
,\,...,\alpha _n$ are different from zero, then
\begin{equation}\label{eq2}
\frac{{\alpha _1  + \alpha _2  + ... + \alpha _\nu  }} {\nu }
\geqslant \sqrt[\nu ]{{\alpha _1  \cdot \alpha _2 ... \cdot \alpha
_\nu  }} \geqslant \frac{\nu } {{\dfrac{1} {{\alpha _1 }} +
\dfrac{1} {{\alpha _2 }} + ... + \dfrac{1} {{\alpha _\nu  }}}}~.
\end{equation}

\noindent Equality is true only when $ \alpha _1  = \alpha _2  = ... =
\alpha _\nu $.
\end{proposition}

For a nice proof, we propose the book of Rassias~\cite{Ra}.\\

\noindent{\footnotesize\textbf{Historical Notation 1\footnote{All Historical Notations are found in the Wikipedia biographies of Cauchy, Schwarz, Bunyakovsky, Holder and Minkowski.
}}. \emph{This
inequality first appeared in Augustin-Louis Cauchy's book
(1789-1857) entitled "Coursd 'Analysis of the Polycole
Polytechnique" (1821) for students at the Ecole Polytechnique in
Paris. The method of proving inequality is a "kind of mathematical
induction". A very beautiful proof of this inequality is presented
in the book (see [Mathematics, Th. Rassias p. 271]) and is due to
the great Hungarian mathematician George Polya (1887-1985). It is
worth noting that this inequality has been given a particularly
large amount of evidence.}}

\noindent{\footnotesize\textbf{Historical Notation 2}. \emph{Baron
Augustin-Louis Cauchy, French mathematician and pioneer of
Analysis, was born on August 21, 1789 and died on May 23, 1857. He
graduated as a civil engineer in 1807 at the age of 18 from the
Polytechnic School of Paris, where he was a professor. He began to
formulate and prove the theorems of infinite calculus in a
rigorous way, rejecting any heuristic principle of the generality
of algebra used by older writers. He defined the sequence in
infinite terms, gave several important theorems in complex
analysis, and began the study of transpositional groups in
abstract algebra. Cauchy, being a thoughtful mathematician,
greatly influenced modern mathematicians as well as later ones.})}

\smallbreak Before formulating the Cauchy-Schwarz-Bunyakovsky
inequality we will recall the concept of the inner product in
$\mathbb{R}^n ,\,\,n \geqslant 1$, and its basic properties. In
Calculus we had to define the concept of absolute value in order
to study distance in real line. The same can happen in the vector
space $\mathbb{R}^n$. We know that the $n$ vectors $ e_1 = \left(
{1,0,...,0} \right),\,e_2  = \left( {0,1,...,0} \right),\,...,e_n
= \left( {0,0,...,1} \right)$ of $\mathbb{R}^n$ consist a basis of
it and in fact the normal basis (that is, it is unit and
perpendicular to two) and they generalize its usual orthogonal
basis $ B = \left\{ {i,\,j,\,k} \right\}\,,\,\,\,\left( {i = e_1
,\,j = e_2 ,\,k = e_3 } \right) $ of $\mathbb{R}^3$.\\
We know that for two possible vectors $ x = \left( {x_1 ,\,x_2
,x_3 } \right) $
 and $
y = \left( {y_1 ,\,y_2 ,y_3 } \right) $ of $\mathbb{R}^3$
  their  inner product is defined as
  \[
x \cdot y = x_1 y_1  + x_2 y_2  +  x_3 y_3 \in \Bbb{R}.
\]
We can now define in a similar way the inner product of two
vectors of $\mathbb{R}^n$ for any $n$. So if $x = \left( {x_1
,\,...,x_n } \right) \in \mathbb{R}^n$ and $y = \left( {y_1
,\,...,y_n } \right) \in \mathbb{R}^n$
 extending the   definition we gave for the inner product
of two vectors of $\mathbb{R}^3$ for any $n$ we can define the
inner product of the vectors of $\mathbb{R}^n$ in a similar way via
\[
x \cdot y = x_1 y_1  + x_2 y_2  + ... + x_n y_n  = \sum\limits_{i
= 1}^n {x_i y_i }~.
\]
The length or Euclidean norm of $ x = \left( {x_1 ,\,...,x_n }
\right)$ is defined to be
 the non-negative real \[
\left\| x \right\| = \sqrt {x \cdot x}  = \sqrt {x_1^2  + ... +
x_n^2 }~.
\]
The pair $ \left( {\mathbb{R}^n ,\left\| {\, \cdot \,} \right\|}
\right) $
 is called the Euclidean space of dimension $n$. We often use the notation $\langle x,y\rangle$ instead $x\cdot y$
and the notation $|x|$ instead $\|x\|$.

The basic properties of the inner product are described in the
following  proposition.

\begin{proposition} \label{pro2.2} \rm
 For all $ x,y,z\, \in
\mathbb{R}^n$, the following properties hold
\[
\begin{gathered}
   (\rm i) \,\,\langle x,y \rangle =\langle y,x \rangle  \hfill \\
   (\rm {ii})\,\, \langle x + z, y \rangle=   \langle x , y \rangle +\langle z,y \rangle\hfill \\
   (\rm {iii})\,\,\langle x ,x \rangle\geqslant 0 \hfill \\
   (\rm {iv}) \,\,\langle x ,x \rangle =0 \Leftrightarrow x = 0 \Leftrightarrow x_1  = x_2  = ... = x_n  = 0 \hfill \\
\end{gathered}
\]
\end{proposition}

The proof of properties is achieved by simply applying the definition
of the inner product. \smallbreak

\subsection{Cauchy-Bunyakovsky-Schwarz Inequality} The $\rm
CBS$ inequality is expressed as follows.

\begin{theorem} \label{th2.1}\rm
Let $ x,\,y \in \mathbb{R}^n$. Then
\begin{equation}\label{eq3}
\left| {x \cdot y} \right| \leqslant \left\| x \right\| \cdot
\left\| y \right\|~.
\end{equation}
The equality in $\rm CBS$ inequality holds if and only if the
vectors $x$ and $y$ are collinear.
\end{theorem} \emph{Proof.}
If $x=0$ or $y=0$  then inequality applies trivially as equality.
Assume that $x\neq0$ and $y\neq0$.

$\rm (i)$ Let's say that $\left\| x \right\|$=1= $\left\| y
\right\|$. Then from the obvious inequality $ \left( {\left| {x_i
} \right| - \left| {y_i } \right|} \right)^2  \geqslant 0 $ we get
that
\[
x_i ^2  + y_i ^2  - 2\left| {x_i y_i } \right| \geqslant 0
\]
from which it follows that
\[|x_i  y_i|\leq\frac{x_i^2+y_i^2}{2},\,i=1,2,...,n.\]
Therefore,
\[
\left| {\sum\limits_{i = 1}^n {x_i y_i } } \right| \leqslant
\sum\limits_{i = 1}^n {\left| {x_i y_i } \right|}  \leqslant
\frac{1} {2}\sum\limits_{i = 1}^n {\left( {x_i^2  + y_i^2 }
\right)}  = \frac{1} {2}\sum\limits_{i = 1}^n {x_i^2 }  + \frac{1}
{2}\sum\limits_{i = 1}^n {y_i^2 }  = \frac{1} {2} + \frac{1} {2} =
1~,
\]
that is
\[
\left| {x \cdot y} \right| = \left| {\sum\limits_{i = 1}^n {x_i
y_i } } \right| \leqslant 1 = 1 \cdot 1 = \left| x \right|
\left| y \right|.
\]
Thus, in this case   inequality applies.

$\rm (ii)$ General case: Let $ x = \left( {x_1 ,\,...,x_n }
\right) \in \mathbb{R}^n $ and $ y = \left( {y_1 ,\,...,y_n }
\right) \in \mathbb{R}^n$
 two  possible vectors. We set \\$$
X_i    = \frac{{x_i }} {{\left\| x \right\|}},\,\,  Y_i  =
\frac{{y_i }} {{\left\| y \right\|}} ,  i=1,2,..,n.$$  It follows
from the case $\rm (i)$  that
 \[
\left| {\sum\limits_{i = 1}^n {X_i  Y_i  } } \right| = \left|
{\sum\limits_{i = 1}^n {\left( {\frac{{x_i }} {{\left\| x
\right\|}} \cdot \frac{{y_i }} {{\left\| y \right\|}}} \right)} }
\right| \leqslant \sum\limits_{i = 1}^n {\frac{{\left| {x_i }
\right|}} {{\left\| x \right\|}} \cdot \frac{{\left| {y_i }
\right|}} {{\left\| y \right\|}}}  \leqslant 1~.
\]
Because
\[
\sum\limits_{i = 1}^n {X_i^2 }  = \sum\limits_{i = 1}^n {\left(
{\frac{{x_i }} {{\left\| x \right\|}}} \right)^2  = \sum\limits_{i
= 1}^n {\frac{{x_i ^2 }} {{\left\| x \right\|^2 }} = } } \frac{1}
{{\left\| x \right\|^2 }}\sum\limits_{i = 1}^n {x_i ^2  = }
\frac{1} {{\left\| x \right\|^2 }}\left\| x \right\|^2  = 1
\]
 and similarly
\[
\sum\limits_{i = 1}^n {Y_i^{'2} }  = 1~.
\]
It follows that
\[
\left| {x \cdot y} \right| = \left| {\sum\limits_{i = 1}^n {x_i
y_i } } \right| \leqslant 1 \cdot 1 = \left\| x \right\|\left\| y
\right\|.
\]
If the vectors $x$ and $y$  are linear, that is, if we have for
someone then we easily find the equality in the $\rm CBS$
inequality.

We now assume that $\left| {x \cdot y} \right| = \left\| x
\right\|\left\| y \right\| $. If one of $x$ and $y$ is $0$, let
$x$ be, then we set $\lambda=0$ and observe that $x=\lambda y$. We
now assume that $ x = \left( {x_1 ,\,...,x_n } \right) \ne 0 $
 and $
y = \left( {y_1 ,\,...,y_n } \right) \ne 0 $. Then the equality $
\left| {x \cdot y} \right| = \left\| x \right\|\left\| y
\right\|\, $ is written as \[ \left| {\sum\limits_{i = 1}^n
{\frac{{\left| {x_i } \right|}} {{\left\| x \right\|}} \cdot
\frac{{\left| {y_i } \right|}} {{\left\| y \right\|}}} } \right| =
1
\]
and posing $ X_i  = \dfrac{{x_i }} {{\left\| x \right\|}} $ and $
Y_i  = \dfrac{{y_i }} {{\left\| y \right\|}} $, yields
\[
\left| {\sum\limits_{i = 1}^n {X_i  \cdot Y_i } } \right| = 1,
\]
where as we proved above $ \sum\limits_{i = 1}^n {X_i^2 }  = 1 =
\sum\limits_{i = 1}^n {Y_i^2 } $.\\ Therefore
\begin{equation}\label{eq4}
1 = \left| {\sum\limits_{i = 1}^n {X_i  \cdot Y_i } } \right|
\leqslant \sum\limits_{i = 1}^n {\left| {X_i  \cdot Y} \right|}
\leqslant \frac{1} {2}\sum\limits_{i = 1}^n {X_i^2 }  + \frac{1}
{2}\sum\limits_{i = 1}^n {Y_i^2 }  = 1.
\end{equation}
We conclude that the numbers $ X_i ,Y_i ,\,i = 1,2,...,n $ are
identical, since the equation holds in the triangular inequality,
that is
\[
\left| {\sum\limits_{i = 1}^n {X_i  \cdot Y_i } } \right| =
\sum\limits_{i = 1}^n {\left| {X_i  \cdot Y_i } \right|}~.
\]
By (\ref{eq4}), also, it   follows that necessary
\begin{equation}\label{eq5}
|X_i  Y_i|  = \frac{1} {2}X_i^2  + \frac{1} {2}Y_i^2 ,\,\,i =
1,\,2,...,n.
\end{equation}

$\rm(a)$  If the common sign of the numbers $ X_i ,Y_i ,\,i =
1,2,...,n $ is positive then from the equalities (\ref{eq5}) we
conclude that
\[
X_i Y_i  = \frac{1} {2}\left( {X_i^2  + Y_i^2 } \right)
\Leftrightarrow \left( {X_i  - Y_i } \right)^2  = 0
\Leftrightarrow X_i  = Y_i ,\,\,\,\forall \,\,i = 1,\,2,...,n\,
\]
or
\[
\frac{{x_i }} {{\left\| x \right\|}} = \frac{{y_i }} {{\left\| y
\right\|}} \Leftrightarrow x_i  = \frac{{\left\| x \right\|}}
{{\left\| y \right\|}} \cdot y_i  \Leftrightarrow x =
\frac{{\left\| x \right\|}} {{\left\| y \right\|}} \cdot y
\]
and
\[
\lambda  = \frac{{\left\| x \right\|}} {{\left\| y \right\|}}.
\]

$\rm(b)$ If the common sign of the numbers $ X_i ,Y_i ,\,i =
1,2,...,n $ is negative then from the equalities (\ref{eq5}) we
find that
\[
-X_i Y_i  = \frac{1} {2}\left( {X_i^2  + Y_i^2 } \right)
\Leftrightarrow \left( {X_i + Y_i } \right)^2  = 0 \Leftrightarrow
X_i  = -Y_i ,\,\,\,\forall \,\,i = 1,\,2,...,n\,
\]
and
\[
\lambda  =- \frac{{\left\| x \right\|}} {{\left\| y \right\|}}.
\]
 \mbox{ }\hfill $\Box$\\

\noindent{\footnotesize\textbf{Historical Notation 3}. \emph{Karl
Hermann Amandus Schwarz, a German mathematician known for his work
on complex analysis, was born on January 25, 1843 in Hermsnstorf,
Silesia (southern Poland) and died on November 30, 1921 in Berlin.
Schwarz initially studied chemistry in Berlin, but Kummer and Karl
Weierstrass encouraged him to switch to mathematics. He received
his Ph.D. from the University of Berlin in 1864 under the
supervision of Kummer and Weierstrass. Between 1867 and 1869, he
worked at the University of Halle and then at the Swiss Federal
Polytechnic. From 1875 he worked at the University of Gottingen
and was involved in complex analysis, differential geometry and
calculus of change. }}

\noindent{\footnotesize\textbf{Historical Notation 4}. \emph{Viktor
Yakovlevich Bunyakovsky, a Russian mathematician known for his
work in theoretical engineering and number theory, was born on
December 16, 1804 in Bar, Russia, and died on November 12, 1889 in
St. Petersburg. He was a member and later vice-president of the
Russian Academy of Sciences. He graduated from the Sorbonne
University in 1824 and received his doctorate under Cauchy. He
proved the $\rm CBS$ inequality for the infinite case in 1859,
long before Schwarz dealt with it.}}

\subsection{Geometrical Interpretation of the $\rm CBS$ Inequality}

From the $\rm CBS$ inequality it follows that for non-zero vectors
$x$ and $y$ in $\Bbb{R}^n$ the quotient \[ \frac{{x \cdot y}}
{{\left\| x \right\|\left\| y \right\|}}
\]
 belongs to the interval $[-1,1]$, i.e.
\[
 - 1 \leqslant \frac{{x \cdot y}}
{{\left\| x \right\|\left\| y \right\|}} \leqslant 1~.
\]
This defines the angle $\theta$ between $x$  and $y$  by the
formula
\begin{equation}\label{eq6}
\cos \theta  = \frac{{x \cdot y}} {{\left\| x \right\| \cdot
\left\| y \right\|}},\,\theta  \in \left[ {0,\pi } \right].
\end{equation}
Moreover, if the vectors $x$ and $y$  belong to $\Bbb{R}^3$ then
$\rm CBS$ inequality  follows from the law of
the cosine of trigonometry.

Applying the law of cosines to the triangle with vertice
$O(0,0,0)$ and sides the vectors $x$ and $y$ it will result that:
\begin{equation}\label{eq7}
\left\| {x - y} \right\|^2  = \left\| x \right\|^2  + \left\| y
\right\|^2  - 2\left\| x \right\|\left\| y \right\|\cos \theta ,
\end{equation}
where $\theta$ the angle of $x$  and $y$. From (\ref{eq7})  using the properties of the inner product we get
successively
\[
\begin{aligned}
~&\left\| {x - y} \right\|^2  = \left\| x \right\|^2  + \left\| y \right\|^2  - 2\left\| x \right\| \cdot \left\| y \right\| \cdot \cos \theta \,   \hfill \\
\Leftrightarrow ~&\left( {x - y} \right) \cdot \left( {x - y} \right) = \left\| x \right\|^2  + \left\| y \right\|^2  - 2\left\| x \right\|\left\| y \right\|\cos \theta \,  \hfill \\
\Leftrightarrow ~&x^2  - 2x \cdot y + y^2  = x^2  + y^2  - 2\left\| x \right\|\left\| y \right\|\cos \theta \,  \hfill \\
\Leftrightarrow ~&- 2x \cdot y =  - 2\left\| x \right\|\left\| y \right\|\cos \theta \,  \hfill \\
\Leftrightarrow ~&x \cdot y = \left\| x \right\|\left\| y \right\|\cos \theta \,. \hfill \\
\end{aligned}
\]
Therefore,
\[
\left| {\cos \theta } \right| = \frac{{\left| {x \cdot y}
\right|}} {{\left\| x \right\|\left\| y \right\|}} \leqslant 1
\Leftrightarrow \left| {x \cdot y} \right| \leqslant \left\| x
\right\|\left\| y \right\|.
\]

If the vectors $x$ and $y$ are collinear, that is $x=\lambda y$,
for some  $\lambda\in\Bbb{R}$, then $\theta=0$ or $\theta=\pi$,
therefore $\cos \theta=1$ and so $ \left| {x \cdot y} \right| =
\left\| x \right\|\left\| y \right\|$. If the vectors $x$ and $y$ are   not collinear, that is, linearly
independent, then $ \theta  \in \left( {0,\pi } \right)$ therefore
$\left| {\cos \theta } \right| < 1$  and so $ \left| {x \cdot y}
\right| < \left\| x \right\|\left\| y \right\|$. So, we come to the part where we talk about the inner  product of
non-zeros vectors
\[
x \cdot y = \left\| x \right\|\left\| y \right\|\cos \theta
,\,\,\,x,y \in \mathbb{R}^n \,,
\]
where $\theta  \in \left[ {0,\pi } \right]$ is the unique angle between $x$ and $y$. Two vectors $x,y\in\Bbb{R}^n\backslash\{0\}$ are called
\emph{orthogonal}  if $x\cdot y=0\Leftrightarrow
\theta=\dfrac{\pi}{2}$.

\subsection{Triangular   Inequality}

We can prove our well-known triangular inequality in $\Bbb{R}^n$
on the basis of the $\rm CBS$ inequality. That is, we will prove
that if $x,y\in \Bbb{R}^n$, then we have
\begin{equation}\label{eq8}
\left\| {x + y} \right\| \leqslant \left\| x \right\| + \left\| y
\right\|.
\end{equation}

Indeed using the $\rm CBS$ inequality (\ref{eq3})
 we have
\begin{eqnarray*}
 \left\| {x + y} \right\|^2  &=& \left( {x + y} \right) \cdot \left( {x + y} \right) \hfill \\
& =& x \cdot x + x \cdot y + y \cdot x + y \cdot y \hfill \\
& = & \left\| x \right\|^2  + 2x \cdot y + \left\| y \right\|^2  \hfill \\
& \leqslant & \left\| x \right\|^2  + 2\left\| x \right\|\left\| y \right\| + \left\| y \right\|^2  \hfill \\
& = & \left( {\left\| x \right\| + \left\| y \right\|} \right)^2 .
\end{eqnarray*}
Taking square roots in each member results in the required
inequality. \smallbreak Based on the $\rm CBS$ inequality we can
get the following very useful proposition.

\begin{proposition}\label{pro2.3} \rm
If $(x_1, ,...,x_n)$ and $(y_1,...,y_n)$   are two sequences of
real numbers, then the following inequality holds
\begin{equation}\label{eq9}
\left( {\sum\limits_{i = 1}^n {x_i^2 } } \right)\left(
{\sum\limits_{i = 1}^n {y_i^2 } } \right) \geqslant \left(
{\sum\limits_{i = 1}^n {x_i y_i } } \right)^2,
\end{equation}
with the equality to be valid if and only if the sequences $(x_1,
,...,x_n)$ and $(y_1,...,y_n)$ are proportional term by term, that
is, if there exists  a fixed number $\lambda$ such that $ x_k   =
\lambda y_\kappa  \,\forall \,k  \in \left\{ {1,2,...,n,...}
\right\} $.
 \end{proposition}

  $\rm CBS$ inequality is known to be one of the most basic and
useful inequalities in algebra, and because of its particular
weight it has greatly preoccupied mathematicians. This has
resulted in a very large number of proofs of great interest. Then
we will present some of the ones we have chosen based on the
diversity, the brevity and the beauty that they
present and of course these criteria are subjective.

In the Appendix we present some of the most beautiful and useful
proofs of the $\rm CBS$ Inequality.

  \smallbreak $\rm CBS$ inequality,  as mentioned and above,
without any  doubt,  is one of the most widely used and important
inequalities in all mathematics. If, for example, we write it  in
the  form
\[
x_1 y_1  + x_2 y_2  + ... + x_n y_n  \leqslant \sqrt {x_1^2  +
x_2^2  + ... + x_n^2 } \sqrt {y_1^2  + y_2^2  + ... + y_n^2 } ,
\]
one of the most important and immediate conclusions that emerges
from this is that from $ \sum\limits_{i = 1}^\infty  {x_i^2 }  <
\infty$  and $ \sum\limits_{i = 1}^\infty  {y_i^2 }  < \infty $
implies that $ \sum\limits_{i = 1}^\infty  {\left| {x_i y_i }
\right|}  < \infty $.

\smallbreak We can generalize the above CBS inequality in
arbitrary linear spaces.   Let for example  $V$ be an arbitrary
real vector space. If we can define  a function in $V\times V$
which is defined by the map $ ( x,y ) \mapsto
 \langle  x,y  \rangle$ in a such way that the $\langle  x,y  \rangle$ satisfies the properties
 in Proposition \ref{pro2.2} then the $  \langle  x,y  \rangle$ is an inner product and $ (
{V, \langle { \cdot , \cdot }  \rangle }  ) $ is a real space with
inner product.

If we go one step further and think that the integration of a
function in a space is nothing more than a sum in the continuous
case, the question arises whether the $\rm CBS$ inequality holds
in this continuous case. That is, if we replace the sums with
integrals  arises the following  inequality
\[
\int_a^b {f\left( x \right)g\left( x \right)dx}  \leqslant \left(
{\int_a^b {f^2 \left( x \right)dx} } \right)^{\frac{1} {2}} \left(
{\int_a^b {g^2 \left( x \right)dx} } \right)^{\frac{1} {2}}.
\]
The question therefore is whether such inequality applies. The
answer to this question is affirmative  and in fact we can get a
more general inequality the so-called H\"older inequality. In
other words, the following proposition is valid.

\begin{proposition}\label{pro2.4} \rm (\emph{H\"older's inequality})
Let $p, q$ be two real numbers   such that $1<p<+\infty$ and
$\dfrac{1}{p}+\dfrac{1}{q}=1$. Let, also, the integrable functions
$f,g\,:\left[ {\alpha ,\,\beta } \right] \to \left[ {0,\, + \infty
} \right] $. Then the following inequality  holds
\[
\int_a^b {f\left( x \right)g\left( x \right)dx}  \leqslant \left(
{\int_a^b {f^p \left( x \right)dx} } \right)^{\frac{1} {p}} \left(
{\int_a^b {g^q \left( x \right)dx} } \right)^{\frac{1} {q}},
\]
A special case of the latter is   taken for $ p = q = \dfrac{1} {2}
$. The positive exponents $p$ and $q$ are called \emph{conjugate}
exponents if $ \dfrac{1} {p} + \dfrac{1} {q} = 1 $ or equivalent $p
+ q = pq $.
 \end{proposition}
 \noindent{\footnotesize\textbf{Historical Notation 5.} \emph{We
mention that this inequality first appeared in print in Victor
Yacovlevich Bunyakovsky's Memoire which was published by the
Imperial Academy of Sciences of St. Petersburg in 1859.}}\\

In aim to prove H\"older's inequality we will assume
that
\[ 0 < \left( {\int_a^b {f^p \left( x \right)dx} }
\right)^{\frac{1} {p}} <  + \infty, \,\,\, 0 < \left( {\int_a^b
{g^q \left( x \right)dx} } \right)^{\frac{1} {q}}  <  + \infty,
\]
 and that
 \[
\left( {\int_a^b {f^p \left( x \right)dx} } \right)^{\frac{1} {p}}
\ne 0,\,\,\,\left( {\int_a^b {g^q \left( x \right)dx} }
\right)^{\frac{1} {q}}  \ne 0,
\]
 that is, that the functions $f^p$ and $g^q$ are integrable with non-zero
integrals because these cases are trivial.

For the proof we will
need the following proposition.

\begin{proposition}\label{pro2.5} \rm
Let $p, q$ be two real numbers   such that $1<p<+\infty$ and
$\dfrac{1}{p}+\dfrac{1}{q}=1$. Then, for all $x,y\geq0$, the
following inequality holds
\[
xy \leqslant \frac{{x^p }} {p} + \frac{{y^q }} {q}.
\]
 \end{proposition}
\emph{Proof.} We observe that if $x=0$
or $y=0$ the inequality
holds, then we assume that $x>0$ and $y>0$. Dividing inequality by members with $y^q$ we obtain the inequality \[
\frac{x} {{y^{q - 1} }} \leqslant \frac{{x^p }} {{py^q }} +
\frac{1} {q}
\]
and if we put in the latter $ t = \dfrac{{x^p }} {{y^q }} $,
because of $ \dfrac{1} {p} + \dfrac{1} {q} = 1$, we have $(q-1)p=q$
and so it turns out that \[ \frac{x} {{y^{q - 1} }} = \left(
{\frac{{x^p }} {{y^{\left( {q - 1} \right)p} }}} \right)^{\frac{1}
{p}}  = \left( {\frac{{x^p }} {{y^{pq - p} }}} \right)^{\frac{1}
{p}}  = \left( {\frac{{x^p }} {{y^q }}} \right)^{\frac{1} {p}}  =
t^{\frac{1} {p}}~,
\]
so we get the inequality \[ t^{\frac{1} {p}}  \leqslant \frac{1}
{p}t + \frac{1} {q},\,\,t > 0.
\]
The function $ f\left( t \right) = t^{\frac{1} {p}}  - \frac{1}
{p}t - \frac{1} {q},\,\,t > 0 $
 is differentiable with $
f'\left( t \right) = \frac{1} {p}\left( {t^{ - \frac{1} {q}}  - 1}
\right) $.\\
We observe that $f'(1)=0$, $ f'\left( t \right) < 0,\,\,\,\forall
\,\,t > 1 $ and $ f'\left( t \right) > 0,\,\,\,\forall \,\,t \in
(0,1) $ so  $ f\left( t \right) \leqslant f(1)\, = 0\,\,\forall
\,\,t > 0 $,   and then \[ t^{\frac{1} {p}}  - \frac{1} {p}t -
\frac{1} {q} \leqslant 0,\,\,\forall \,\,t > 0,
\]
where the equality holds  only when $t=1$.
\mbox{ }\hfill $\Box$

\smallbreak\noindent\emph{Proof of H\"older's inequality}. We set
$ A = \left( {\int_a^b {f^p \left( x \right)dx} }
\right)^{\frac{1} {p}} $ and $ B = \left( {\int_a^b {g^q \left( x
\right)dx} } \right)^{\frac{1} {q}} $.  So from the above
proposition for the numbers $  f ( x  )/A  $
 and $
\ g ( x  ) /B  $
 we get
 \[
\frac{{f\left( x \right)}} {A}\frac{{g\left( x \right)}} {B}
\leqslant \frac{ (f(x)/A )^p  } {p} + \frac{(g(x)/B  )^q } {q} =
\frac{{\left( {f\left( x \right)} \right)^p }} {{pA^p }} +
\frac{{\left( {g\left( x \right)} \right)^q }} {{qB^q }}~.
\]
By  integration by parts  from the last we obtain \[
\begin{gathered}
  \,\frac{1}
{{AB}}\int\limits_\alpha ^\beta  {f\left( x \right)} g\left( x
\right)dx \leqslant \frac{1} {{pA^p }}\int\limits_\alpha ^\beta
{\left( {f\left( x \right)} \right)^p dx}  + \frac{1}
{{qB^q }}\int\limits_\alpha ^\beta  {\left( {g\left( x \right)} \right)^q dx}  \hfill \\
  \,\,\,\,\,\,\,\,\,\,\,\,\,\,\,\,\,\,\,\,\,\,\,\,\,\,\,\,\,\,\,\,\,\,\,\,\,\,\,\,\,\,\,\,\,\,\, = \frac{1}
{{pA^p }}A^p  + \frac{1} {{qB^q }}B^q    = \frac{1} {p} + \frac{1}
{q} = 1, \hfill \\
\end{gathered}
\]
so \[ \int\limits_\alpha ^\beta  {f\left( x \right)} g\left( x
\right)dx \leqslant AB \Leftrightarrow \,\int\limits_\alpha ^\beta
{f\left( x \right)} g\left( x \right)dx \leqslant \left(
{\int\limits_\alpha ^\beta  {f^p \left( x \right)} }
\right)^{\frac{1} {p}} \left( {\int\limits_\alpha ^\beta  {g^q
\left( x \right)} } \right)^{\frac{1} {q}}.
\]
If we use the classic notation for the norm in $L^p$ spaces, that
is, if we set $$ \left\| f \right\|_p  = \left(
{\int\limits_\alpha ^\beta  {f^p \left( x \right)} }
\right)^{\frac{1} {p}} ,\,\,\,1 \leqslant p <  + \infty ,$$ the
last inequality takes the form \[ \left\| {fg} \right\|_1  =
\left\| f \right\|_p \left\| g \right\|_q.
\] \mbox{ }\hfill $\Box$\\

\noindent{\footnotesize \textbf{Historical Notation 6}. \emph{
Ludwig Otto Holder, German mathematician, was born on December 22,
1859 in Stuttgart and died on August 29, 1937. He studied at the
present University of Stuttgart where he was a student of
Kronecker, Weirstrass and Kummer. He received his doctorate from
the University of Tubingen in 1882. He is known for many theorems
such as the Holder Inequality, the Jordan-Holder theorem and the
abnormal external automorphism of the symmetric group.}}\\

 Completing
this generalization in the continuous case, that is, in spaces, we
can get the generalization of the triangular inequality called
Minkowski inequality.

\begin{proposition}\label{pro2.6} \rm
Let $1<p<+\infty$ and the integrable functions $f,\,g\,:\,\left[
{\alpha ,\,\beta } \right] \to \left[ {0,\, + \infty } \right] $.
Then,
\[
\left( {\int_a^b {\left( {f + g} \right)^p \left( x \right)dx} }
\right)^{\frac{1} {p}}  \leqslant \left( {\int_a^b {f^p \left( x
\right)dx} } \right)^{\frac{1} {p}}  + \left( {\int_a^b {g^p
\left( x \right)dx} } \right)^{\frac{1} {p}}
\]
\emph{Proof.} We notice that
\begin{eqnarray*}
  \left[ {f\left( x \right) + g\left( x \right)} \right]^p  &\leqslant &
  \left[ {2\max \left\{ {f\left( x \right),g\left( x \right)} \right\}} \right]^p \hfill \\
 & = & 2^p \max \left\{ {f^p \left( x \right),g^p \left( x \right)} \right\} \hfill \\
 & \leqslant & 2^p \left[ {f^p \left( x \right) + g^p \left( x \right)} \right]
\end{eqnarray*}
or
\[
\left( {f\left( x \right) + g\left( x \right)} \right)^p \leqslant
2^p f^p \left( x \right) + 2^p g^p \left( x \right),
\]
 so with integration by parts we get this inequality \[
\int_a^b {\left( {f + g} \right)^p \left( x \right)dx}  \leqslant
2^p \int_a^b {f^p \left( x \right)dx}  + 2^p \int_a^b {g^p \left(
x \right)dx}.
\]
We can assume that $ \int_a^b {f^p \left( x \right)dx}  <  +
\infty $, $ \int_a^b {g^p \left( x \right)dx}  \leqslant  + \infty
$, and  then due to the last inequality arises $ \int_a^b {\left(
{f + g} \right)^p \left( x \right)dx}  <  + \infty $
 and so in this case the inequality we want to prove is obvious.
 
 Let now be the conjugate exponent $q$ of $p$. Then, $ \left( {p -
1} \right)q = p$,
 so from H\"older's inequality
we get  \[ \int_a^b {f\left( {f + g} \right)^{p - 1} \left( x
\right)dx} \leqslant \left( {\int_a^b {f^p \left( x \right)dx} }
\right)^{\frac{1} {p}} \left( {\int_a^b {\left( {f + g} \right)^p
\left( x \right)dx} } \right)^{\frac{1} {q}}
\]and \[
\int_a^b {g\left( {f + g} \right)^{p - 1} \left( x \right)dx}
\leqslant \left( {\int_a^b {g^p \left( x \right)dx} }
\right)^{\frac{1} {p}} \left( {\int_a^b {\left( {f + g} \right)^p
\left( x \right)dx} } \right)^{\frac{1} {q}}.
\]
Adding by members the last two inequalities, we reach \[ \int_a^b
{\left( {f + g} \right)^p \left( x \right)dx}  \leqslant \left[
{\left( {\int_a^b {f^p \left( x \right)dx} } \right)^{\frac{1}
{p}}  + \left( {\int_a^b {g^p \left( x \right)dx} }
\right)^{\frac{1} {p}} } \right]\left( {\int_a^b {\left( {f + g}
\right)^p \left( x \right)dx} } \right)^{\frac{1} {q}}
\]
and from the latter, by division by members, we obtain the inequality \[
\left( {\int_a^b {\left( {f + g} \right)^p \left( x \right)dx} }
\right)^{1 - \frac{1} {q}}  \leqslant \left( {\int_a^b {f^p \left(
x \right)dx} } \right)^{\frac{1} {p}}  + \left( {\int_a^b {g^p
\left( x \right)dx} } \right)^{\frac{1} {p}}
\]
and because $ 1 - \dfrac{1} {q} = \dfrac{1} {p} $ we finally arrive to \[
\left( {\int_a^b {\left( {f + g} \right)^p \left( x \right)dx} }
\right)^{\frac{1} {p}}  \leqslant \left( {\int_a^b {f^p \left( x
\right)dx} } \right)^{\frac{1} {p}}  + \left( {\int_a^b {g^p
\left( x \right)dx} } \right)^{\frac{1} {p}}.
\]
The last inequality, again using the classic notation of
$L^p$-norm, can be written in the form \[ \left\| {f + g}
\right\|_p  \leqslant \left\| f \right\|_p  + \left\| g
\right\|_p,
\]
which is the generalization of triangular inequality in spaces.
\mbox{ }\hfill $\Box$
 \end{proposition}

\noindent{\footnotesize \textbf{Historical Notation 7}.
\emph{Herman Minkowski, German mathematician, was born on 22 June
1864 in Lithuania and died on 12 January 1909 in Gottingen. He
received his PhD from the University of Konigsbgerg at the age of
21. He taught at the Federal Institute of Technology (ETH) in
Zurich and in the last years of his life at the University of
Gottingen, where he had the mathematics course created for him by
David Hilbert. But Schwarz needed a two-dimensional form of Cauchy
inequality. In particular, he had to show that if $S \subset
\mathbb{R}^2$ and $f,g:S \to \mathbb{R}$, then the double
integrals $A = \iint_S {f^2 dxdy}$, $B = \iint_S {fg dxdy}$, and
$A = \iint_S {g^2 dxdy}$ must satisfy the inequality $ \left| B
\right| \leqslant \sqrt A  \cdot \sqrt C $
 and he had to know that the inequality is strict unless the
functions are proportional. Schwarz would face a problem for a
number of reasons, including the fact that the severity of the
inequality can be lost using integrals. So he had to find a
different solution, which he did as he discovered a proof that
went down in history. He based his proof on an impressive
observation. He observed  that the real polynomial $ p\left( t
\right) = \iint_S {\left( {t\,f\left( {x,y} \right) + g\left(
{x,y} \right)} \right)}^2 dxdy = At^2  + 2Bt + C
$
 is always non-negative and that $p(t)$  it is strictly positive unless $f$ and
and $g$ are proportional. The binomial formula tells us that the
coefficients must satisfy $B^2  \leqslant AC$
 and unless $f$ and $g$ are proportional then
we have strict inequality $B^2 < AC$. So from a simple algebraic
observation Schwarz discovered everything he needed to know.}}

\section {A Cauchy-Bunyakovsky-Schwarz  Type  Inequality}

In this section we present a Cauchy-Schwarz type inequality, which
arose from the study of mechanical properties of elastic rod (see
in \cite{Bi-Bi}). Firstly, we present the discreet case and
secondly its integral version.

\begin{theorem}\label{th3.1} \rm
Let $x=(x_1,...,x_n)$ and $y=(y_1,...,y_n)$ be arbitrary elements
of $\mathbb{R}^n$. Then the following   holds

\begin{eqnarray}\label{eq10}
\nonumber\left( {\sum\limits_{i = 1}^n {p_i } } \right)\left(
{\sum\limits_{i = 1}^n {p_i x_i } y_i } \right) - \left(
{\sum\limits_{i = 1}^n {p_i x_i } } \right)\left( {\sum\limits_{i
= 1}^n {p_i } y_i } \right) \\
\leqslant  \sqrt {\left(
{\sum\limits_{i = 1}^n {p_i } } \right)\left( {\sum\limits_{i =
1}^n {p_i x_i^2 } } \right) - \left( {\sum\limits_{i = 1}^n {p_i
x_i } } \right)^2 } \\ 
\nonumber\times\sqrt {\left(
{\sum\limits_{i = 1}^n {p_i } } \right)\left( {\sum\limits_{i =
1}^n {p_i y_i^2 } } \right) - \left( {\sum\limits_{i = 1}^n {p_i
y_i } } \right)^2 },
\end{eqnarray}
for any positive real numbers $p_1, p_2,...,p_n$ (weights).  The
above relation (\ref{eq10}) becomes an equality if and only if
there exist $a,b\in \mathbb{R}$ such  that $ax+by=c$, a constant
vector, i.e. if there exists a linear combination of the vectors
$x$ and $y$ which is a constant vector.
\end{theorem}
\emph{Sketch of the proof}. We present the basic steps of the
proof. (A proof in detail and some important remarks are available
in \cite{Bi-Bi}).

We define  the function of two variables $\langle . ,
.\rangle:\mathbb{R}^n\times\mathbb{R}^n\to \mathbb{R}$ defined by
\begin{eqnarray}\label{eq11}
\left\langle {x,y} \right\rangle  = \left( {\sum\limits_{i = 1}^n
{p_i } } \right)\left( {\sum\limits_{i = 1}^n {p_i x_i } y_i }
\right) - \left( {\sum\limits_{i = 1}^n {p_i x_i } } \right)\left(
{\sum\limits_{i = 1}^n {p_i } y_i } \right)~.
\end{eqnarray}
It is easy to prove that the function $\langle . ,
.\rangle:\mathbb{R}^n\times\mathbb{R}^n\to \mathbb{R}$ is a
bilinear form, i.e. it satisfies the  properties of Proposition \ref{pro2.2}. These properties can be easy  verified by taking into account the definition (\ref{eq11}).
Concerning the first property it can be proved using the CBS
inequality (\ref{eq3}), that is
\begin{eqnarray*}\label{eq17.21}
\left\langle {x,x} \right\rangle & = & \left( {\sum\limits_{i =
1}^n {p_i } }
  \right)\left( {\sum\limits_{i = 1}^n {p_i x_i^2 } _i } \right) - \left(
   {\sum\limits_{i = 1}^n {p_i x_i } } \right)\left( {\sum\limits_{i = 1}^n {p_i } x_i } \right) \hfill \\
 & = &\left[ {\sum\limits_{i = 1}^n {\left( {\sqrt {p_i }
   } \right)^2 } } \right]\left[ {\sum\limits_{i = 1}^n {\left( {\sqrt {p_i } x_i }
    \right)^2 } } \right] - \left( {\sum\limits_{i = 1}^n {\sqrt {p_i } \sqrt {p_i x_i } } } \right)^2 \, \geqslant
    0.
\end{eqnarray*}
The inequality becomes equality,  i.e.~$\left\langle {x,x}
\right\rangle=0$,  if and  only if the vectors $(\sqrt{p_1},
\sqrt{p_2},..., \sqrt{p_n})$  and  $(\sqrt{p_1}x_1, \sqrt{p_2}x_2,
... , \sqrt{p_n}x_n)$ are linear dependent, that is,
$x=(x_1,x_2,...,x_n)$ is a constant vector. With these
considerations inequality (\ref{eq11}) takes  the form
\[
\left\langle {x,y} \right\rangle   \leqslant \sqrt{\left\langle
{x,x} \right\rangle} \sqrt{\left\langle {y,y} \right\rangle}
,\,\,\,\forall \,x,y \in \mathbb{R}^n~.
\]
which is  inequality (\ref{eq2}).

The integral version of the CBS type inequality is presented by
the following proposition.  (We will relax  the notation
$\int_\Omega f(x)dx$ and  simply  write  $\int_\Omega f$).

\begin{proposition}\label{pro3.2} \rm
Let the  functions $f,g,p: \Omega\subset\mathbb{R}^n\to\mathbb{R}$
be continuous on a compact domain $\Omega$, and let $p$ be
positive. The following inequality holds
\begin{eqnarray}\label{eq12}
\hspace{-5pt} \int\limits_\Omega  p \int\limits_\Omega  {pfg}  -\hspace{-3pt}
\int\limits_\Omega  p f\int\limits_\Omega  p g \leqslant \sqrt
{\int\limits_\Omega  p \int\limits_\Omega  p f^2  -\hspace{-3pt} \left(
{\int\limits_\Omega  p f} \right)^2 } \hspace{-3pt}\sqrt {\int\limits_\Omega
p \int\limits_\Omega  p g^2  -\hspace{-3pt} \left( {\int\limits_\Omega  p g}
\right)^2},
\end{eqnarray}
where the inequality becomes equality if and only  if there exist
$\alpha,\beta\in \mathbb{R}$ such that $\alpha f + \beta g=c$, a
constant.
 \end{proposition}
\textbf{Notation}.  We can imagine that inequality (\ref{eq12}) is
obtained directly from (\ref{eq10}) by replacing the discrete
``sums" by continuous ones, i.e. by integrals, and what remains to
be shown is that the steps in the proof of the continuous case are
exactly the same as those of the discrete one. \\
\emph{Sketch of the proof}. We set
\begin{eqnarray}\label{eq13}
\left\langle {f,g} \right\rangle  = \int\limits_\Omega  p
\int\limits_\Omega  p fg - \int\limits_\Omega  p f -
\int\limits_\Omega  p g
\end{eqnarray}
It easy to verify  that the properties $\rm(i)-(v)$ in Proposition
\ref{pro2.2} are satisfied   and  inequality  (\ref{eq12}) can be
written
\begin{eqnarray}\label{eq14}
\left\langle {f,g} \right\rangle  \leqslant \sqrt {\left\langle
{f,f} \right\rangle } \sqrt {\left\langle {g,g} \right\rangle
},\,\,\,\forall \,f,g \in C(\Omega),
\end{eqnarray}
which leads  to the proof of  the result.

\section {Strengthened Cauchy-Bunyakovsky-Schwarz
Inequality for Elasticity}

In this section we present the strengthened CBS inequality,
\begin{equation}\label{eq17.21}
|\langle u,v\rangle|\leq \gamma \sqrt{\langle
u,u\rangle}\sqrt{\langle v,v\rangle},\,\,u\in U,\,\,v\in V,\,\,
U\cap V=\{0\},
\end{equation}
where $U$ and $V$ are two linear subspaces $\langle.,.\rangle$ is
the bilinear form corresponding to the variational formulation of
the problem and $\gamma$ is a constant between 0 and 1 depending
only on the spaces $U$ and $V$. (In fact the constant $\gamma$
could take the values 0 and 1, but if $\gamma=0$ then the vectors
$u$ and $v$ are orthogonal and for $\gamma=1$ arises   the
classical CBS inequality). The estimate of the constant that
appears in the second term of this inequality and in particular
its optimal value (that is, the smallest value it can taken to
hold for all $u\in U$ and $v\in V$) is of fundamental practical
importance in solving the problems that appear and especially in
solving 2D and 3D elasticity problems, which are the ones that
interest us. For example, Margenov in \cite{Mar} gave estimates on
the 2D   elasticity problem on a triangular mesh and proved for a
uniform mesh of right isosceles triangles, that the constant
$\gamma^2$ is bounded above by $3/4$ uniformly on the mesh.
Achchab and Maitre in \cite{Ac-Ma} proved that this result remains
true for every triangular mesh. In the case of a 3D problem, Jung
and Maire in \cite{Ju-Ma} reported that numerical experiments on a
particular triangulation and specifically in the refinement of a
standard tetrahedron showed that the constant $\gamma^2$ is
bounded by $9/10$. Achchab et al. in \cite{Ach-Axe} generalized
this result and proved that it remains valid for an arbitrary
tetrahedron mesh. In addition, Achchab et al. proves that the
constant $\gamma$ plays significant  role in posteriori error
estimates for elasticity problems.

In this section we present  the basic theorem needed for
the establishment of the strengthened CBS inequality, which is needed for the estimation of the
constant $\gamma$. (For more details and for the proofs see also
the work of Achchab et al.\cite{Ach-Axe})

\begin{theorem}\label{th4.1} \rm
Given a finite-dimensional Hilbert space $H$, an inner product
$\langle,.,\rangle$ on it and the subspaces $U$ and $V$ of $H$
with $U\cap V=\{0\}$, there exist a constant $\gamma\in (0,1)$
depending only  on $U$ and $V$ such that for all $u\in U$ and
$v\in V$ the following strengthened CBS inequality holds
\begin{equation}\label{eq17.21}
|\langle u,v\rangle|\leq\gamma \|u\| \|v\|,
\end{equation}
where the norm induced by the inner product.
\end{theorem}

\begin{theorem}\label{th4.1} \rm
Let $A$ be a symmetric, positive semidefine $2\times 2$ block
matrix \[ A = \left( {\begin{array}{*{20}c}
   {A_{11} } & {A_{12} }  \\
   {A_{21} } & {{\text{A}}_{{\text{22}}} }  \\
 \end{array} } \right),
\]
where $A_{11}$ is invertible and where the partitioning components
to $U$ , $V$, which are the spaces of vectors with only non-zero
first, and second components $u\in U$, $v\in V$ of the form $ u =
\left( {\begin{array}{*{20}c}
   {u_1 }  \\
   0  \\
 \end{array} } \right)
$ and
$
v = \left( {\begin{array}{*{20}c}
   0  \\
   {v_2 }  \\
 \end{array} } \right)
,$ respectively. Assuming that the kernel of $A$, N(A), is
included in $V$, the optimal value of the constant $\gamma$ is
given by
\[
\gamma ^2  = \mathop {\mathop {\sup }\limits_{u \in U} }\limits_{v
\in V\backslash N\left( A \right)} \frac{{\left( {u^T Av}
\right)^2 }} {{\left( {u^T Au} \right)\left( {v^T Av} \right)}}~.
\]
\end{theorem}

Let $\Omega$ be a bounded open connected subset of $\Bbb R^3$ with
a Lipschitz continuous boundary $\Gamma$. Let, also, $\Gamma_0$ a
measurable subset of $\Gamma$ and $\Gamma_1$ its complement on
$\Gamma$. We consider the following 3-dimensional boundary value
problem of linear elasticity
\[
\begin{gathered}
  \sum\limits_{j = 1}^3 {\frac{{\partial \sigma _{ij} \left( u \right)}}
{{\partial x_j }}}  + f_i \,\,\,\,\,\mathrm{in}\,\,\,\Omega ,\,\,\,i = 1,2,3 \hfill \\
  u_i  = 0\,\,\,\,\,\,\,\,\,\,\,\,\,\,\,\,\,\,\,\,\,\,\,\,\,\,\,\,\,\,\,\mathrm{on}\,\,\Gamma _0 ,\,\,i = 1,2,3 \hfill \\
  \sum\limits_{j = 1}^3 {\sigma _{ij} \left( u \right)v_j }  = g_i \,\,\,\,\mathrm{on}\,\,\Gamma _1 ,\,\,i = 1,2,3. \hfill \\
\end{gathered}
\]
$v=(v_1,v_2,v_3)$ is the external normal to $\Omega$, with the
classical constitutive law for isotropic material with L\'ame
modula $\lambda$ and $\mu$
\[
c_{ij} \left( u \right) = \lambda \left( {\sum\limits_{k = 1}^3
{\varepsilon _{kk} } \left( u \right)} \right)\delta _{ij}  + 2\mu
\varepsilon _{ij} \left( u \right),
\]
where
\[
\varepsilon _{ij} \left( u \right) = \frac{1} {2}\left(
{\frac{{\partial u_i }} {{\partial x_j }} + \frac{{\partial u_j }}
{{\partial x_i }}} \right).
\]
We introduce the space \[ V\left( \Omega  \right) = \left\{ {v \in
\left[ {H^1 \left( \Omega  \right)} \right]^3 ;\,\,v =
0\,\,\mathrm{on}\,\,\Gamma _0 } \right\}
\]
For $ f \in \left[ {L^2 \left( \Omega  \right)} \right]^3 $ and $
\,g \in \left[ {L^2 \left( {\Gamma _0 } \right)} \right]^3 $ the
variational formulation of the problem is given by
$$ \mathrm{Find}\,\,u\in
V(\Omega)\,\,\mathrm{such}\,
\mathrm{that}\,\,\,a(u,v)=L(v),\,\,\forall\,v\in \Omega,$$ where
\[
a\left( {u,v} \right) = \lambda \int_\Omega  {div\left( u \right)}
div\left( v \right)dx + 2\mu \sum\limits_{i,j = 1}^3 {\varepsilon
_{ij} \left( u \right)\varepsilon _{ij} \left( v \right)} dx
\]
and
\[
L\left( v \right) = \int_\Omega  {fv} dx + \int_{\Gamma _1 } {gv}
d\sigma.
\]
Finally, we define the two forms $a_1(,.,)$ and $a_2(,.,)$ by
$$a_1(u,v)=\int_\Omega \rm{div} (u) \rm{div}(v)dx$$
and $$ \;\;\qquad a_2 \left( {u,v} \right) = \sum\limits_{i,j =
1}^3 {\int_\Omega {c_{ij} } } \left( u \right)c_{ij} \left( v
\right)dx .$$ Under the above considerations the following theorem
holds.
\begin{theorem}\label{th0.1} \rm
The constants $\gamma_1$, $\gamma_2$ of the strengthened CBS
inequality associated the forms $a_1(u,v)$  and $ a_2 \left( {u,v}
\right)  $ satisfies
$$\gamma_i^2\leq\frac{9}{10},\,\,i=1,2.$$
\end{theorem}

%%%%%%%%%%%%%%%%%%%%%%%%%%%%%%%%%%%%%%%
%\begin{appendix}[]
\section*{Appendix}
%\section{Appendices}\index{appendix}
%\smallbreak
%\appendix{\textbf{Appendix A}} \smallbreak\noindent
\emph{Proof 1.}  Consider the double sum $ \sum\limits_{i = 1}^n
{\sum\limits_{j = 1}^n } (  x_i y_j  - x_j y_i  )  ^2 .$
 By deleting the parentheses and grouping the same terms we get
successively
\begin{eqnarray*}
\sum\limits_{i = 1}^n {\sum\limits_{j = 1}^n {\left( {x_i y_j  -
x_j y_i } \right)}  ^2} & = & \sum\limits_{i = 1}^n
{\sum\limits_{j = 1}^n {\left( {x_i ^2 y_j ^2
+ x_j ^2 y_i ^2  - 2x_i x_j y_i y_j } \right)} }  \hfill \\
& =& \sum\limits_{i = 1}^n {\left( {\sum\limits_{j = 1}^n {x_i ^2 y_j ^2 }  + \sum\limits_{j = 1}^n {x_j ^2 y_i ^2 }  - \sum\limits_{j = 1}^n {2x_i x_j y_i y_j } } \right)}  \hfill \\
& =&\sum\limits_{i = 1}^n {\left( {x_i ^2 \sum\limits_{j = 1}^n {y_j ^2 }  + y_i ^2 \sum\limits_{j = 1}^n {x_j ^2 }  - 2x_i y_i \sum\limits_{j = 1}^n {x_j y_j } } \right)}  \hfill \\
& = &\sum\limits_{i = 1}^n {\left( {x_i ^2 \sum\limits_{j = 1}^n {y_j ^2 } } \right)}  + \sum\limits_{i = 1}^n {\left( {y_i ^2 \sum\limits_{j = 1}^n {x_j ^2 } } \right)\, - 2\sum\limits_{i = 1}^n {\left( {x_i y_i \sum\limits_{j = 1}^n {x_j y_j } } \right)} }  \hfill \\
& = &\sum\limits_{j = 1}^n {y_j ^2 } \sum\limits_{i = 1}^n {x_i ^2 }  + \sum\limits_{j = 1}^n {x_j ^2 } \sum\limits_{i = 1}^n {y_i ^2 \, - 2\sum\limits_{j = 1}^n {x_j y_j } \sum\limits_{i = 1}^n {x_i y_i } }  \hfill \\
&= &\sum\limits_{i = 1}^n {x_i^2 \sum\limits_{j = 1}^n {y_j^2 } }  + \sum\limits_{i = 1}^n {y_i^2 \sum\limits_{j = 1}^n {x_j^2 } }  - 2\sum\limits_{j = 1}^n {x_i y_i } \sum\limits_{j = 1}^n {x_j y_j }  \hfill \\
 & = &
   \,\,2\left( {\sum\limits_{i = 1}^n {x_i^2 } } \right)\left( {\sum\limits_{i = 1}^n {y_i^2 } } \right) - 2\left(
   {\sum\limits_{j = 1}^n {x_i y_i } } \right)^2
\end{eqnarray*}
Because the left-hand side of the equation is a sum of squares of
real numbers is greater or equal  to zero, so \[ \left(
{\sum\limits_{i = 1}^n {x_i^2 } } \right)\left( {\sum\limits_{i =
1}^n {y_i^2 } } \right) \geqslant \left( {\sum\limits_{i = 1}^n
{x_i y_i } } \right)^2,
\]
that is, the inequality (\ref{eq9}). \mbox{ }\hfill $\Box$\\

\smallbreak\noindent\emph{Proof 2.} Consider the next trinomial
with respect to $t$
\[
f\left( t \right) = \left( {\sum\limits_{i = 1}^n {x_i^2 } }
\right)t^2  - 2\left( {\sum\limits_{i = 1}^n {x_i y_i } } \right)t
+ \sum\limits_{i = 1}^n {y_i^2 }.
\]
Doing the operations in the second member we notice that
\begin{eqnarray*}
f\left( t \right) &=& \left( {\sum\limits_{i = 1}^n {x_i^2 } }
 \right)t^2  - 2\left( {\sum\limits_{i = 1}^n {x_i y_i } }
 \right)t + \sum\limits_{i = 1}^n {y_i^2 }  \hfill \\
&= &\sum\limits_{i = 1}^n {x_i^2
  } t^2  - 2\sum\limits_{i = 1}^n {x_i y_i } t +
  \sum\limits_{i = 1}^n {y_i^2 }  \hfill \\
& = &\sum\limits_{i = 1}^n
  {\left( {x_i t} \right)^2 }  - 2\sum\limits_{i =
   1}^n {\left( {x_i t} \right)y_i }  + \sum\limits_
   {i = 1}^n {y_i^2 }  \hfill \\
& = &\sum\limits_{i = 1}^n
   {\left[ {\left( {x_i t} \right)^2  - 2\left( {x_i t} \right)y_i  + y_i^2 } \right]}  \hfill \\
& = &\sum\limits_{i = 1}^n
  {\left( {x_i t - y_i } \right)^2 }  \geqslant 0.
\end{eqnarray*}
Since $f\left( t \right) \geqslant 0$
 for every $t\in \Bbb{R}$ the  discriminant  of $f(t)$  is negative so
it holds that \[ \left( {\sum\limits_{i = 1}^n {x_i y_i } }
\right)^2  - \left( {\sum\limits_{i = 1}^n {x_i^2 } }
\right)\left( {\sum\limits_{i = 1}^n {y_i^2 } } \right) \leqslant
0 \Leftrightarrow \left( {\sum\limits_{i = 1}^n {x_i y_i } }
\right)^2  \leqslant \left( {\sum\limits_{i = 1}^n {x_i^2 } }
\right)\left( {\sum\limits_{i = 1}^n {y_i^2 } } \right),
\]
that is, the inequality (\ref{eq9}) holds. \mbox{ }\hfill
$\Box$\\

\smallbreak\noindent\emph{Proof 3.} We assume that the $\rm CBS$
inequality  (\ref{eq4}) is true for every $ \left( {x_1 ,\,...,x_n
} \right) $ and $ \left( {y_1 ,\,...,y_n } \right)$. Based on this
hypothesis we will come to the true. If $ \sum\limits_{i = 1}^n
{x_i^2 }  = 0 $ or $ \sum\limits_{i = 1}^n {y_i^2 }  = 0 $, then $
x_1  = x_2  = ... = x_n  = 0$    or $ y_1  = y_2  = ... = y_n  =
0$  and then the conclusion is obvious. Now we assume that $ X_n
= \sum\limits_{i = 1}^n {x_i^2 }  \ne 0 $ and $ Y_n  =
\sum\limits_{i = 1}^n {y_i^2 }  \ne 0 $ and we set  $ a_i  =
\dfrac{{x_i }} {{\sqrt {X_n } }} $ and $ \,b_i  =
\dfrac{{y_i }} {{\sqrt {Y_n } }}$, for each $i=1,2,...,n$. Then we will have
\[
 \sum\limits_{i = 1}^n {a_i^2 }  =
\sum\limits_{i = 1}^n {\left( {\frac{{x_i }} {{\sqrt {X_n } }}}
\right)^2 }  = \sum\limits_{i = 1}^n {\frac{{x_i ^2 }} {{X_n }} =
} \frac{1} {{X_n }}\sum\limits_{i = 1}^n {x_i ^2  = } \frac{1}
{{X_n }}X_n  = 1
\]
and \[ \sum\limits_{i = 1}^n {b_i ^2 }  = 1.
\]
Since we have assumed that the initial inequality (\ref{eq4}) is
true for all $n$-groups we apply it to $ \left( {a_1 ,\,...,a_n }
\right) $ and $ \left( {b_1 ,\,...,b_n } \right) $ and  then we
get \[ 1 \cdot 1 \geqslant \left( {\sum\limits_{i = 1}^n {a_i b_i
} } \right)^2
\]
from which inequality arises
\[
 a_1 b_1  + a_2 b_2  + ... + a_nb_n  \leqslant 1.
\]
From the latter we take successively
\begin{eqnarray*}
2\left( {a_1 b_1  + a_2 b_2  + ... + a_n b_n } \right) &\leqslant& 2\,\,\, \Leftrightarrow  \hfill \\
2\left( {a_1 b_1  + a_2 b_2  + ... + a_n b_n } \right) &\leqslant& 1 + 1 \Leftrightarrow  \hfill \\
2\left( {a_1 b_1  + a_2 b_2  + ... + a_n b_n } \right) &\leqslant&
\left( {a_1 ^2  + a_2 ^2  + ... + a_n ^2 } \right) + \left( {b_1
^2  + b_2 ^2  + ... + b_n ^2 } \right),
\end{eqnarray*}
thus \[\left( {a_1  - b_1 } \right)^2  + \left( {a_2  - b_2 }
\right)^2
    + ... +
  \left( {a_n  - b_n } \right)^2   \geqslant  0,
\]
which is obviously true. \mbox{ }\hfill $\Box$\\

\smallbreak\noindent\emph{Proof 4.} Let \[ X = \sqrt {x_1^2  +
x_2^2 + ... + x_n^2 } ,\,\,\,Y = \sqrt {y_1^2  + y_2^2  + ... +
y_n^2 }
\]
 From the inequality of the arithmetic-geometric mean,
we have
\begin{eqnarray*}
 \sum\limits_{i = 1}^n {\frac{{x_i y_i }}
{{XY}}}  &\leqslant& \frac{1} {2}\sum\limits_{i = 1}^n {\left(
{\frac{{x_i^2 }} {{X^2 }} + \frac{{y_i^2 }}
{{Y^2 }}} \right)}  \hfill \\
& =& \frac{1} {2}\left( {\sum\limits_{i = 1}^n {\frac{{x_i^2 }}
{{X^2 }}}  + \sum\limits_{i = 1}^n {\frac{{y_i^2 }}
{{Y^2 }}} } \right) \hfill \\
 & = & \frac{1}
{2}\left( {\frac{1} {{X^2 }}\sum\limits_{i = 1}^n {x_i^2 }  +
\frac{1}
{{Y^2 }}\sum\limits_{i = 1}^n {y_i^2 } } \right) \hfill \\
&=& 1,
\end{eqnarray*}
so it will be \[ \sum\limits_{i = 1}^n {x_i y_i }  \leqslant {\rm
X}\Upsilon  = \sqrt {x_1^2  + x_2^2  + ... + x_n^2 } \sqrt {y_1^2
+ y_2^2  + ... + y_n^2 }.
\]
Therefore
\[
\left( {\sum\limits_{i = 1}^n {x_i y_i } } \right)^2  \leqslant
\left( {\sum\limits_{i = 1}^n {x_i^2 } } \right)\left(
{\sum\limits_{i = 1}^n {y_i^2 } } \right).
\]
\mbox{ }\hfill $\Box$\\

\smallbreak\noindent\emph{Proof 5.} Let

%\null\hfill  \hspace*{\fill}
%\begin{eqnarray*}
%					&A_n  &= x_1^2  + x_2^2  + ... + x_n^2, \hfill \\ 
%					&B_n  &= y_1^2  + y_2^2 + ... + y_n^2 \hfill \\ 
%\noindent\mathrm{and} 	&C_n  &= x_1 y_1
%+ x_2 y_2 + ... + x_n y_n.
%\end{eqnarray*}

\begin{tabular}[\textwidth]{@{}l@{\hspace{100pt}}c@{~}c@{~}l}
								&$A_n$	&=& $x_1^2  + x_2^2  + ... + x_n^2$, \\ 
								&$B_n$	&=& $y_1^2  + y_2^2 + ... + y_n^2$  \\ 
\hspace{-18pt} and \hspace*{\fill} \null\hfill	& 		&	&						\\ 
\hspace*{\fill} \null\hfill	&$C_n$  &=& $x_1 y_1+ x_2 y_2 + ... + x_n y_n$.
\end{tabular}

It follows from the inequality of the arithmetic-geometric mean
that
\begin{eqnarray*}
  \frac{{\sum\limits_{i = 1}^n {x_i^2 } \sum\limits_{i = 1}^n {y_i^2 } }}
{{\left( {\sum\limits_{i = 1}^n {x_i y_i } } \right)^2 }} + 1 &=&
\frac{{A_n B_n }} {{C_n^2 }} + 1 = \frac{{\sum\limits_{i = 1}^n
{x_i^2 } B_n }} {{C_n^2 }} + \frac{{\sum\limits_{i = 1}^n {y_i^2 }
}} {{B_n }} = \sum\limits_{i = 1}^n {\left( {\frac{{x_i^2 B_n }}
{{C_n^2 }} + \frac{{y_i^2 }}
{{B_n }}} \right)} a \hfill \\
 & \geqslant & \sum\limits_{i = 1}^n {\left( {2\frac{{x_i^2 B_n }}
{{C_n^2 }} \cdot \frac{{y_i^2 }} {{B_n }}} \right)}  = 2\frac{1}
{{C_n }}\sum\limits_{i = 1}^n {x_i y_i }  = 2\frac{1}
{{C_n }}C_n  = 2~. \hfill \\
\end{eqnarray*}
Therefore
\[
A_n B_n  \geqslant C_n^2,
\]
that is
\[
\left( {x_1^2  + x_2^2  + ... + x_n^2 } \right)\left( {y_1^2  +
y_2^2  + ... + y_n^2 } \right) \geqslant \left( {x_1 y_1  + x_2
y_2  + ... + x_n y_n } \right)^2.
\]
 \mbox{ }\hfill $\Box$\\

\smallbreak\noindent\emph{Proof 6.} We will prove the $\rm CBS$
inequality (\ref{eq9}) by mathematical induction. Starting the
induction for 1, the hypothesis is trivial. Note that 
\begin{flalign}
\left(
{x_1 y_1  + x_2 y_2 } \right)^2  &= x_1^2 y_1^2  + 2x_1 y_1 x_2 y_2 + x_2^2 y_2^2  \nonumber\\
&\leqslant x_1^2 y_1^2  + x_1^2 y_2^2  + x_2^2 y_1^2 + x_2^2 y_2^2   \nonumber\\
&= \left( {x_1^2  + x_2^2 } \right)\left( {y_1^2  + y_2^2 } \right), \nonumber
\end{flalign}
which shows that the   inequality  (\ref{eq9}) holds for $n=2$.
Assume that the   inequality (\ref{eq9}) holds for $n=k$, that is,
\[
\left( {\sum\limits_{i = 1}^k {x_i y_i } } \right)^2  \leqslant
\left( {\sum\limits_{i = 1}^k {x_i^2 } } \right)\left(
{\sum\limits_{i = 1}^k {y_i^2 } } \right)
\]
 and we will prove that it is true for $n=k+1$. We have
\begin{eqnarray*}
\sqrt {\sum\limits_{i = 1}^{k + 1} {x_i^2 } } \cdot \sqrt
{\sum\limits_{i = 1}^{k + 1} {y_i^2 } } & = &
\sqrt {\sum\limits_{i = 1}^k {x_i^2  + x_{k + 1}^2 } }  \cdot \sqrt {\sum\limits_{i = 1}^k {y_i^2  + y_{k + 1}^2 } }  \hfill \\
& \geqslant & \sqrt {\sum\limits_{i = 1}^k {x_i^2 } }  \cdot \sqrt {\sum\limits_{i = 1}^k {y_i^2 } }  + \left| {x_{k + 1} y_{k + 1} } \right| \hfill \\
 & \geqslant & \sum\limits_{i = 1}^k {\left| {x_i y_i } \right|}  + \left| {x_{k + 1} y_{k + 1} } \right| = \sum\limits_{i = 1}^{k + 1} {\left| {x_i y_i } \right|} . \hfill \\
\end{eqnarray*}
 This means that the   inequality (\ref{eq9}) holds for $n=k+1$, and so we
conclude that $\rm CBS$ holds for all natural numbers. This
completes the proof.
 \mbox{ }\hfill $\Box$

\smallbreak\noindent\emph{Proof 7.} Let
\begin{eqnarray*}
X &=& \left\{ {x_1 y_1 ,...,x_1 y_n ,x_2 y_1 ,...,x_2 y_n ,...,x_n y_1 ,...,x_n y_n } \right\}, \hfill \\
  Y &=& \left\{ {x_1 y_1 ,...,x_1 y_n ,x_2 y_1 ,...,x_2 y_n ,...,x_n y_1 ,...,x_n y_n } \right\}, \hfill \\
  Z &=& \left\{ {x_1 y_1 ,...,x_1 y_n ,x_2 y_1 ,...,x_2 y_n ,...,x_n y_1 ,...,x_n y_n } \right\}, \hfill \\
  W &=& \left\{ {x_1 y_1 ,...,x_n y_1 ,x_1 y_2 ,...,x_n y_2 ,...,x_1 y_n ,...,x_n y_n } \right\}. \hfill \\
\end{eqnarray*}
It is easy to observe that the $X$ and $Y$ have the same
classification, while the $Z$ and $W$ have a mixed classification.
Applying the rearrangement  inequality, we have
\[
\begin{gathered}
  \left( {x_1 y_1 } \right)\left( {x_1 y_1 } \right) + ... + \left( {x_1 y_n } \right)\left( {x_1 y_n } \right) + \left( {x_2 y_1 } \right)\left( {x_2 y_1 } \right) + ... + \left( {x_2 y_n } \right)\left( {x_2 y_n } \right) + ... \hfill \\
   + \left( {x_n y_1 } \right)\left( {x_n y_1 } \right) + ... + \left( {x_n y_n } \right)\left( {x_n y_n } \right)   \hfill \\
 \geqslant \left( {x_1 y_1 } \right)\left( {x_1 y_1 } \right) + ... + \left( {x_1 y_n } \right)\left( {x_n y_1 } \right) + \left( {x_2 y_1 } \right)\left( {x_1 y_2 } \right) + ... + \left( {x_2 y_n } \right)\left( {x_n y_2 } \right) + ... \hfill \\
   + \left( {x_n y_1 } \right)\left( {x_n y_1 } \right) + ... + \left( {x_n y_n } \right)\left( {x_n y_n } \right), \hfill \\
\end{gathered}
\]
which can be simplified and written as follows \[ \left( {x_1^2  +
x_2^2  + ... + x_n^2 } \right)\left( {y_1^2  + y_2^2  + ... +
y_n^2 } \right) \geqslant \left( {x_1 y_1  + x_2 y_2  + ... + x_n
y_n } \right)^2
\]  The last
inequality is the desired one. \mbox{ }\hfill $\Box$\\

\smallbreak\noindent\emph{Proof 8.} Because if $ x = \left(
{0,0,...,0} \right) $  or $ y= \left( {0,0,...,0} \right) $
 inequality (\ref{eq4}) is always valid, we
consider the non-zero vectors $ x = \left( {x_1 ,x_2 ,...,x_n }
\right)$ and $y = \left( {y_1 ,y_2 ,...,y_n } \right) $.
 Then for any real number $t$, the
following holds \[ \left| {x + ty} \right|^2  = \left( {x + ty}
\right) \cdot \left( {x + ty} \right) = x \cdot x + 2\left( {x
\cdot y} \right)t + \left( {y \cdot y} \right)t^2  = \left| x
\right|^2  + 2\left( {x \cdot y} \right)t + \left| y \right|^2
t^2,
\]
so inequality \[ \left| x \right|^2  + 2\left( {x \cdot y}
\right)t + \left| y \right|^2 t^2  \geqslant 0
\]
holds  for every  $ t\, \in \,\mathbb{R}$. But the first term of
the last inequality is a trinomial with respect to $t$ and because
it is greater than or equal to $0$ for each $t\in \Bbb{R}$ and its
coefficient is $ \left| y \right|^2  > 0$ (because we have assumed
that $ y = \left( {y_1 ,y_2 ,...,y_n } \right) \ne 0 $ its
discriminant  will be less than or equal to zero. Therefore, \[
4\left( {x \cdot y} \right)^2  - 4\left| x \right|^2 \left| y
\right|^2  \leqslant 0
\]
 so
it will be
\[
\left( {x \cdot y} \right)^2  \leqslant \left| x \right|^2 \left|
y \right|^2.
\]
Using the notations \[ x \cdot y = x_1 y_1  + x_2 y_2  + ... + x_n
y_n ,\,\,\left| x \right|^2  = \sum\limits_{i = 1}^k {x_i^2 }
,\,\,\left| y \right|^2  = \sum\limits_{i = 1}^k {y_i^2 } ,
\]
we result the following
\[
\left( {\sum\limits_{i = 1}^n {x_i y_i } } \right)^2  \leqslant
\left( {\sum\limits_{i = 1}^n {x_i^2 } } \right)\left(
{\sum\limits_{i = 1}^n {y_i^2 } } \right),
\]
which is  the inequality (\ref{eq9}). \mbox{ }\hfill $\Box$\\

\smallbreak\noindent\emph{Proof 9.} Consider the vectors \[ x =
\left( {x_1 ,x_2 ,...,x_n } \right),\,\,y = \left( {y_1 ,y_2
,...,y_n } \right).
\]
 From the formula of inner product $x \cdot y = \left| x \right|\left| y \right|\cos \left( {x,y}
\right)$
 we conclude that $
x \cdot y \leqslant \left| x \right|\left| y \right| $. Using the
notations \[ x \cdot y = x_1 y_1  + x_2 y_2  + ... + x_n y_n
,\,\,\left| x \right|^2  = \sum\limits_{i = 1}^k {x_i^2 }
,\,\,\left| y \right|^2  = \sum\limits_{i = 1}^k {y_i^2 },
\]
the $\rm CBS$ inequality arises. \mbox{ }\hfill
%\end{appendix}\\

\bibliographystyle{unsrtnat}

\bibliography{Cauchy_Inequality_biblio.bib}


%  \bibitem{kour2014real}
%  George Kour and Raid Saabne.
%  \newblock Real-time segmentation of on-line handwritten arabic script.
%  \newblock In {\em Frontiers in Handwriting Recognition (ICFHR), 2014 14th
%        International Conference on}, pages 417--422. IEEE, 2014.

%  \bibitem{kour2014fast}
%  George Kour and Raid Saabne.
%  \newblock Fast classification of handwritten on-line arabic characters.
%  \newblock In {\em Soft Computing and Pattern Recognition (SoCPaR), 2014 6th
%        International Conference of}, pages 312--318. IEEE, 2014.

%  \bibitem{hadash2018estimate}
%  Guy Hadash, Einat Kermany, Boaz Carmeli, Ofer Lavi, George Kour, and Alon
%  Jacovi.
%  \newblock Estimate and replace: A novel approach to integrating deep neural
%  networks with existing applications.
%  \newblock {\em arXiv preprint arXiv:1804.09028}, 2018.

% \end{thebibliography}

\end{document}